\documentclass{tMPH2e}

\usepackage{multirow}
\usepackage{bm}
\usepackage[justification=centering]{caption}
\usepackage{amsmath}
\usepackage{float}
\usepackage{graphicx}

\begin{document}

\doi{10.1080/0026897YYxxxxxxxx}
 \issn{1362–3028}
\issnp{0026–8976}
\jvol{00}
\jnum{00} \jyear{2009} 
\articletype{ORIGINAL ARTICLE}

\title{Generalizing The Mean Spherical Approximation as a Multiscale, Nonlinear Boundary Condition at the Solute--Solvent Interface}

\author{Amirhossein Molavi Tabrizi, Matthew G. Knepley, Jaydeep P. Bardhan}

\maketitle

\begin{abstract}
In this paper we extend the familiar continuum electrostatic model with a perturbation to the usual macroscopic boundary condition. The perturbation is based on the mean spherical approximation (MSA), to derive a multiscale hydration-shell boundary condition (HSBC).  We show that the HSBC/MSA model reproduces MSA predictions for Born ions in a variety of polar solvents, including both protic and aprotic solvents.  Importantly, the HSBC/MSA model predicts not only solvation free energies accurately but also solvation entropies, which standard continuum electrostatic models fail to predict.   The HSBC/MSA model depends only on the normal electric field at the dielectric boundary, similar to our recent development of an HSBC model for charge-sign hydration asymmetry, and the reformulation of the MSA as a boundary condition enables its straightforward application to complex molecules such as proteins.
\bigskip

\bigskip

\end{abstract}


\section{Introduction}

The thermodynamics of solute--solvent interactions play critical roles that range from fundamental chemistry and biology, to nanotechnology and environmental science.  Many applications, particularly in chemical engineering, require understanding the relative properties of a wide range of solvents for a particular system, e.g. at a planar electrode.  In such problems, the breadth of solvents, and also solvent conditions such as temperature and pressure, necessitates simple yet robust theories. One of these theories that provide key physical insights while retaining as much simplicity as possible is the Mean Spherical Approximation (MSA)~\cite{Blum75,Blum78,Stell75,Blum80,BlumVericat80,Tang03}. In contrast, in many biological applications the central challenge is to understand small changes in the solute (e.g. changes in protein conformation), while the solvent itself remains unchanged.  Biological solutions are for the most part dilute aqueous electrolytes, composed of a monovalent salt concentration of approximately 145 mM.  This dilute nature motivated theoreticians to model the electrostatics of biomolecules using the Poisson-Boltzmann equation outside the molecule~\cite{Kirkwood34,Tanford57} and macroscopic dielectric theory for the biomolecule and the molecule--solvent boundary conditions.

The advent of high-performance computer simulations have enabled divergent tracks for studying the solvation of complex biomolecules.  On one hand, large-scale molecular dynamics simulations in explicit solvent offer fully atomistic detail.  On the other hand, efficient, parallelized algorithms for the Poisson-Boltzmann (PB) equation enable large-scale studies with millions of protein atoms treated in atomic detail, while the solvent is treated using continuum theory such as the PB equation~\cite{Sharp90,Baker01,Yokota11}.  However, standard PB models neglect important physical phenomena such as correlations between solvent molecules (including both water and ions)~\cite{Ren12_review,Bardhan12_review}.  Numerous modified PB models for strong electrolytes address ion packing effects and correlations~\cite{Borukhov97,Vlachy99,Grochowski07,Harris14_Boschitsch_Fenley,Gillespie15} but fewer address the equivalent challenges for water itself.  Popular approaches for water's non-bulk response include modeling saturation via nonlinear dielectric theory~\cite{Sandberg02_first_LD,Gong09,Hu12_Wei_nonlinear_Poisson_BJ,Abrashkin07} and solvent dipole correlations via nonlocal dielectric theory~\cite{Dogonadze74,Hildebrandt04,Fedorov07,Bardhan11_pka}.   Models for water steric effects have been incorporated in continuum models via, inter alia, Ornstein-Zernike integral equations~\cite{Chandler72} and density functional theories~\cite{Ramirez02,Knepley10,Borgis12}, which represent much more sophisticated approaches than classical PB.  Among the important failings of macroscopic continuum models is their failure to reproduce important solvent thermodynamics, e.g. in electron transfer~\cite{Matyushov96,Vath99}, while models such as the MSA succeed~\cite{Vath99}.

Unfortunately, the majority of studies based on the MSA address geometries that are much simpler than atomistic models of proteins or DNA. In this paper, we utilized the MSA model for Born ions to derive a multiscale \textit{hydration-shell boundary condition} (HSBC) for the solute-solvent interface. This HSBC can be used for any non-spherical molecular shapes. This modified boundary condition accounts for solvent molecular-size effects, and appears as a nonlinear perturbation to the usual (macroscopic) dielectric boundary condition. Our approach leads to a continuum model that accurately predicts the solvation free energies and entropies of Born ions, in water as well as in other polar solvents. We call this model the HSBC/MSA because it captures hydration-shell phenomena using the MSA. The new HSBC/MSA continuum model is readily solved for large, complex solutes such as proteins and colloids, and easily implemented in existing software for continuum electrostatics (e.g. DelPhi~\cite{Rocchia01} and APBS~\cite{Baker01}). Other HSBCs based on improved liquid state theories are possible through a similar approach. 

It should be noted that the HSBC/MSA model presented in this paper is not an approximation or simplification of the MSA model, but rather a \textit{localization} of the MSA model. MSA considers the dipole-ion interactions by introducing a correction to the ion radius that is based on the solvent molecular radius. Because the MSA model accounts for the finite size of the solvent's molecules by modifying the radius, Implementing it for non-spherical solutes is not an easy task.The main contribution of this work is to model the effect of a hydration shell as a multiscale nonlinear boundary condition, rather than a function of the radii of solvent molecules. Because computation of the HSBC/MSA model is completely local, it can be used to compute electrostatic interactions among large biomolecules with geometrically complex boundary with the boundary element method. The main computational use of MSA is to determine the single fitting parameter in our model.
This kind of localization has been used inconjunction with MSA to define electrostatic correlation structure in
classical density functional theory~\cite{Gillespie02,Gillespie15}, as well as dielectric correlations~\cite{Gillespie11b}.

The following section introduces fundamental theoretical results regarding the continuum electrostatic model, boundary-integral (surface-charge) interpretations of the Born equation, and calculating Born ion solvation free energies via the MSA. Section~\ref{sec:msa-as-nlbc} contains the central result of this paper: a derivation of a new multiscale HSBC that augments the traditional macroscopic boundary condition with a correction based on the MSA. Section~\ref{sec:results} presents computational results demonstrating the new model's accuracy for monoatomic ions in various solvents, and Section~\ref{sec:discussion} concludes the paper by highlighting opportunities for further refinement.


\section{Theory}\label{sec:theory}

\subsection{A Boundary Integral Equation for the Continuum-Electrostatic Model}

We first introduce the standard macroscopic continuum model and the corresponding boundary integral equation. We assume that the solute is a linear continuum medium with $N_q$ point charges inside, and that the electric potential here, denoted $\phi_{in}$, obeys the Poisson equation
\begin{align}
\nabla^2 \phi_{in}(\bm{r})=-\frac{\rho(\bm{r})}{\epsilon_{in}},
\label{eq:poisson}
\end{align}
where $\bm{r}$ is a point in space, $\epsilon_{in}$ is the solute dielectric constant, and the charge distribution $\rho(\bm{r})=\sum_{i=1}^{N_q} q_i\delta({\bm{r}-\bm{r}_i})$, where $\bm{r}_i$ is the position vector of the $i^{th}$ charge and $\delta$ is the Dirac delta function. The solute-solvent boundary is sharp and denoted by $S$, and we model the exterior solvent as a linear dielectric medium with $\epsilon_{out}\gg\epsilon_{in}$. The electrostatic potential here satisfies the Laplace equation
\begin{align}
\nabla^2 \phi_{out}(\bm{r})=0.
\label{eq:laplace}
\end{align}
Macroscopic dielectric theory and Gauss's law lead to the standard Maxwell boundary conditions (SMBC)
\begin{align}
\phi_{in}(\bm{r}_{_{S}})&=\phi_{out}(\bm{r}_{_{S}});\label{re:bd1}\\
\epsilon_{in}\frac{\partial \phi_{in}}{\partial n}(\bm{r}_{_{S}})&=\epsilon_{out}\frac{\partial \phi_{out}}{\partial n}(\bm{r}_{_{S}}).\label{eq:bd2}
\end{align}
In Eq.~\ref{eq:bd2}, $\partial/\partial n$ is the (outward) normal derivative and we assume that $\phi_{out}\rightarrow 0$ quickly enough as $|\bm{r}|\rightarrow \infty$. This mixed-dielectric Poisson problem is well posed and we can write a boundary integral equation (BIE) for the induced surface charge on $S$, $\sigma(\bm{r})+\bar{\epsilon}\int\sigma(\bm{r}^{\prime})\frac{\partial G}{\partial n}d\bm{r}^{\prime}=-\bar{\epsilon}\sum q_i \frac{\partial G}{\partial n}$, where $\bar{\epsilon}=(\epsilon_{in}-\epsilon_{out})/\left({\frac{1}{2}\left(\epsilon_{in}+\epsilon_{out}\right)}\right)$, and $G$ is the three-dimensional free-space Green's function, $G=1/(4\pi |{\bm{r}-\bm{r}^{\prime}}|)$. Once found, we can then find the potential inside the solute or solvent as desired. The BIE problem is equivalent to the problem of a homogeneous medium with relative permittivity $\epsilon_{in}$ and a surface charge distribution $\sigma(\bm{r})$ on the boundary. the potential in the equivalent problem is
\begin{align}
\phi_{in}(\bm{r})=\sum{q_i}\frac{G(\bm{r};\bm{r_i})}{\epsilon_{in}}+\int_{S}\sigma({\bm{r}^{\prime}})\frac{G(\bm{r};\bm{r}^{\prime})}{\epsilon_{in}}d\bm{r}^{\prime}.
\label{eq:poten-equiv}
\end{align}
Using Eq.~\ref{eq:poten-equiv} we can find the normal derivative of the potential
\begin{align}
\frac{\partial \phi_{in}(\bm{r})}{\partial n}=\frac{1}{\epsilon_{in}}\sum{q_i}\frac{\partial G(\bm{r};\bm{r_i})}{\partial n}+\frac{1}{\epsilon_{in}}\int_{S}\sigma({\bm{r}^{\prime}})\frac{\partial G(\bm{r};\bm{r_i})}{\partial n}d\bm{r}^{\prime},
\label{eq:poten-equiv-der}
\end{align}
and Gauss's law in a homogeneous medium means that the jump condition at the surface is
\begin{align}
\frac{\sigma(\bm{r})}{\epsilon_{in}}=\frac{\partial \phi_{in}(\bm{r})}{\partial n(\bm{r})}-\frac{\partial \phi_{out}(\bm{r})}{\partial n(\bm{r})}.
\label{eq:jump-condition}
\end{align}
Approaching the field point $\bm{r}$ to the surface, we obtain the boundary integral equation for induced charge $\sigma(\bm{r})$ on the surface
\begin{align}
\left( I +\hat{\epsilon}\left( -\frac{1}{2}I+K\right)\right)\sigma=\hat{\epsilon}\sum_{i=1}^{N_q} q_i\frac{\partial G}{\partial n},
\label{eq:jump}
\end{align}
where $\hat{\epsilon}=(\epsilon_{out}-\epsilon_{in})/\epsilon_{out}$, $K$ is the normal electric-field operator and $I$ is the identity. The presented equations are valid for an arbitrarily shaped surface, $S$. In the next sections we discuss the Born and MSA models for ion solvation, explicitly considering spherical ions.


\subsection{A Surface-Charge View of the Born Model for Ion Solvation}

The well-known Born equation expresses the Gibbs solvation free energy for a spherical ion as
\begin{align}
\Delta G^{Born}=\frac{-N_L(z_ie_0)^2}{8\pi \epsilon_0}\left( \frac{1}{\epsilon_{in}}-\frac{1}{\epsilon_{out}}\right) \frac{1}{R},
\label{eq:gibbs-born-1}
\end{align}
where $N_L$ is Avogadro's constant, $z_i$ is the ion valence, $e_0$ the fundamental unit of charge, $\epsilon_0$ the permittivity of free space, $R$ the ion radius, and $\epsilon_{in}$ and $\epsilon_{out}$ the dielectric constants for the ion and solvent. In this work, we assume $\epsilon_{in}=1$. Assuming linear response and writing the solvation free energy in terms of the dielectric boundary charge $\sigma$ from the previous section, we have
\begin{align}
\Delta G=\frac{1}{2}q\phi^{Reac}=\frac{1}{2}q\int_{S}{\frac{\sigma({\bm{r}^{\prime}})}{4\pi |{\bm{r}-\bm{r}^{\prime}}|}}d\bm{r^\prime}=\frac{1}{2}qR\sigma,
\label{eq:gibbs-reac}
\end{align}
where $\phi^{Reac}$ is the reaction potential and $q$ is the charge. The last equality holds because $\sigma({\bm{r}^{\prime}})=\sigma$ due to symmetry. Using Eqs.~\ref{eq:jump} and \ref{eq:gibbs-reac}, and noting that a constant $\sigma(\bm{r}^{\prime})$ is an eigenfunction of $K$ with eigenvalue $1/2$, we have that $(-\frac{1}{2}I+K)\sigma=0$ and so we recover Eq.~\ref{eq:gibbs-born-1} with the surface charge
\begin{align}
\sigma^{Born}=C ~\hat{\epsilon} ~ \frac{-q}{4\pi R^2}=C~ \hat{\epsilon} ~ q~ \frac{\partial{G}}{\partial{n}},
\label{eq:sigma-born}
\end{align}
where $C=N_L/(\epsilon_0 \epsilon_{in})$ and $z_ie_0=q$.


In this classical continuum Born model, the effect of temperature is usually underestimated~\cite{Vath99}. This underestimation occurs because the Gibbs free energy is a function of temperature through $\hat{\epsilon}$ only, and $\epsilon_{out}$ is always much larger than $\epsilon_{in}$.


\subsection{MSA Reference Model}
Macroscopic dielectric theory assumes that the constituent dipoles are infinitesimally small compared to the system.  However, for real solutions, the surrounding solvent molecules are not infinitely small, and the electric field induced by the solute charges disturbs the solvent structure. One of the models that considers the effect of dipole-ion interactions is the MSA \cite{Blum87,Fawcett04}, and in the MSA model the Gibbs solvation free energy of a Born ion is
\begin{align}
\Delta G^{MSA}=\frac{-N_L(z_ie_0)^2}{8\pi \epsilon_0}\left(\frac{1}{\epsilon_{in}}-\frac{1}{\epsilon_{out}}\right)\frac{1}{R+\delta_s},
\label{eq:gibbs-MSA-1}
\end{align}
where $\delta_s$ \textit{depends on the solvent}.  In particular, $\delta_s=R_s/\lambda_s$ where $R_s$ is the radius of the solvent and $\lambda_s$ is the MSA solvent polarization parameter. Importantly, $\lambda_s$ is a function of the bulk permittivity $\epsilon_{out}$, and can be calculated by the Wertheim relationship \cite{Fawcett92, Fawcett04}
\begin{align}
\lambda_s^2(1+\lambda_s)^4=16\epsilon_{out}.
\label{eq:Wertheim}
\end{align}
Similarly to the Born model, we may derive the MSA's effective surface charge using Eq.~\ref{eq:gibbs-MSA-1} and Eq.~\ref{eq:gibbs-reac}:
\begin{align}
\sigma^{MSA}=-\frac{C}{4\pi}~\hat{\epsilon} ~ \frac{q}{R(R+\delta_s)}.
\label{eq:sigma-MSA}
\end{align}
The parameter $\delta_s$ is often assumed to modify the ion radius.  However, in this work we assume that in fact the ion radius remains at the original value $R$.  That is, we do not consider $\delta_s$ to involve modification of the ion radius; instead, we focus on reproducing its \textit{effect} (on the reaction potential in the solute) by finding a boundary condition that generates the MSA surface charge $\sigma^{MSA}$, on the surface of the "actual" ion with radius $R$.
\section{Deriving a Multiscale Boundary Condition from the MSA}\label{sec:msa-as-nlbc}

The approximations in the standard continuum-electrostatic model are especially problematic  at solute-solvent boundaries.  In particular, due to the strong electric field and the disturbed solvent structure, the boundary condition Eq.~\ref{eq:bd2} \textit{is not valid near the boundary!} On the other hand, this boundary condition is easy to apply in general geometries such as proteins at atomic resolution.  In contrast, the much more accurate MSA~\cite{Fawcett92} is difficult to apply in inhomogeneous systems without high symmetry. For example, Eq.~\ref{eq:sigma-MSA} depends explicitly on the ion radius and perturbation $\delta_s$.

We sought to eliminate these dependencies by finding a simple modified version of Eq.~\ref{eq:jump} whose solution might approximate the effective MSA charge density (Eq.~\ref{eq:sigma-MSA}). We considered the modification
\begin{align}
\left( I + h(E_n) +\hat{\epsilon}\left( -\frac{1}{2}I+K\right)\right)\sigma=\hat{\epsilon}\sum_{i=1}^{N_q} q_i\frac{\partial G}{\partial n},
\label{eq:modified-bie}
\end{align}
where the new term is $h(E_n)$, which is a function that relies only on the normal electric field at the boundary just inside the solute, and $E_n=\sum_{i}^{Nq}q_i\frac{\partial G}{\partial n}-\bm{K}\sigma$.  In principle, the function $h$ could depend on the tangential field or higher-order derivatives of the potential, but here we restrict the approximation to use only the normal field. The key feature is that the Born problem's spherical symmetry simplifies Eq.~\ref{eq:modified-bie} to
\begin{align}
(1+h(E_n))\sigma^{HSBC/MSA}=C~ \hat{\epsilon} ~q~ \frac{\partial{G}}{\partial{n}},
\label{eq:NLBC-1}
\end{align}
and we can sample $h(E_n)$ using a ``test set'' of Born ions of varying size. Thus, considering the MSA surface charge $\sigma^{MSA}$ as a target, we defined an exact $h(E_n)$ by substituting it in for $\sigma^{HSBC/MSA}$ in Eq.~\ref{eq:NLBC-1} and sought to find a function that might approximate 
\begin{align}
h(E_n)=\frac{\sigma^{Born}}{\sigma^{MSA}}-1,
\label{eq:NLBC-2}
\end{align}
where $\sigma^{Born}$ and $\sigma^{MSA}$ are presented in Eqs.~\ref{eq:sigma-born} and \ref{eq:sigma-MSA} respectively. For the MSA model, we propose
\begin{align} 
  h(E_n)=\alpha \sqrt{|E_n|},
\label{eq:h(En)}
\end{align}
 where $\alpha$ is a solvent- and tempreture-dependent fitting parameter. Finally by introducing $f(E_n)$,
 \begin{align}
f(E_n)=\frac{\epsilon_{in}}{\epsilon_{out}-\epsilon_{in}}-h(E_n),
\label{eq:h}
\end{align}
we see that the modified boundary condition we call HSBC/MSA is 
\begin{align}
f(E_n)\frac{\partial \phi_{in}}{\partial n}(\bm{r}_{_{S}})&=\left(1+f(E_n)\right)\frac{\partial \phi_{out}}{\partial n}(\bm{r}_{_{S}})
\label{eq:nlbc0}\\
\left(\epsilon_1 - \Delta\epsilon h(E_n) \right)\frac{\partial \phi_{in}}{\partial n}&=\left(\epsilon_2 - \Delta\epsilon h(E_n)\right)\frac{\partial \phi_{out}}{\partial n},
\end{align}
where  we have defined $\Delta \epsilon = \epsilon_2-\epsilon_1$.  HSBC/MSA therefore recovers the classical continuum model as $\alpha\rightarrow 0$, and the Born surface charge in the HSBC/MSA model is
\begin{align} 
\sigma^{HSBC/MSA}=C~\hat{\epsilon}~q~\frac{\partial{G}}{\partial{n}}~\frac{1}{1+h(E_n)}.
\label{eq:NLBC-sigma}
\end{align}
One interesting advantage of the proposed HSBC/MSA model is that temperature effects are included naturally, following the same mechanism as in the original MSA.  Fawcett and Blum have tabulated $\delta_s$ and $d \delta_s/d T$ for numerous polar solvents at $25^\circ$C~\cite{Fawcett92}.  To test the model's robustness over a much wider range of temperatures,  however, we parameterized  $\alpha$ in Eq.~\ref{eq:h(En)} in multiple solvents, for which the dielectric constant has been parameterized as a function of temperature.  For instance, in water, for temperatures between 0$^{\circ}$C and 100$^{\circ}$C, the dielectric constant as a function of temperature has been parameterized as \cite{Malmberg56}
\begin{align}
\epsilon_{W}=(-1.410\times10^{-6})T^3+(9.398\times10^{-4})T^2-0.40008T+87.740;
\label{eq:eps-water}
\end{align}
Using this expression for water's dielectric constant and by minimizing the sum of squared difference between
Eq.~\ref{eq:NLBC-2} and Eq.~\ref{eq:h(En)}, we find $\alpha$ at $T=0^\circ$C, $25^\circ$C, $50^\circ$C,
$75^\circ$C, and $100^\circ$C (Table~\ref{table:alpha-water}). Probing the values for $\alpha$ at different temperatures, we can determine
$\alpha$ as a simple function of temperature, with  $\alpha_{W}(T)=0.000594T+0.670476$. Using this equation for $\alpha_{W}$, one can approximate $h(E_n)$ for water to a maximum
relative error of 9\%. It is worth mentioning that this error increases as the radius decreases and temperature increases.
\begin{table}[H]
\begin{center}
\caption{Optimized values of $\alpha$ for water (W) at different temperatures. The graphical representation of $\alpha$ as a linear function of temperature for different solvents are available in Fig.~\ref{fig:alpha}}
\begin{tabular}{ l| c c c c c }
  & $0^{\circ}$C & $25^{\circ}$C & $50^{\circ}$C & $75^{\circ}$C & $100^{\circ}$C  \\
  \hline
$\alpha $ (\AA) & $0.670785$ & $0.685195$ & $0.699823$ & $0.714839$ & $0.730188$ \\
 \end{tabular}
 \label{table:alpha-water}
\end{center}
\end{table}

The $h(E_n)$ at $T=25^\circ$C and $75^\circ$C are shown in Fig.~\ref{fig:func h(En), water}; the solid lines are from Eq.~\ref{eq:NLBC-2} and the points are from Eq.~\ref{eq:h(En)}. It can be observed that as $E_n$ approaches zero (i.e. $R\rightarrow\infty$), $h(E_n)$ approaches zero and therefore the induced surface charge in HSBC/MSA approaches $C\hat{\epsilon}q\frac{\partial{G}}{\partial{n}}$. This is the surface induced charge of a macroscopic Born model and therefore both MSA and HSBC/MSA recover the macro Born results. Also, we have charge symmetry, the function is even, as expected; there are charge-sign asymmetric MSA versions~\cite{Fawcett92} which we intend to investigate in future work.
\begin{figure}[H]
\centering
 \includegraphics[width=25pc]{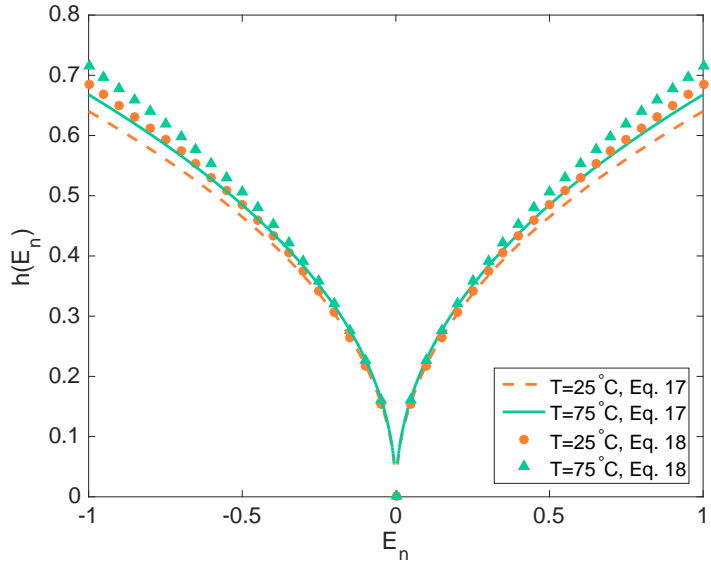}
 \caption{Function $h(E_n)$ for water (W) at $T=25^\circ$C and $T=75^\circ$C. The solid lines are based on Eq.~\ref{eq:NLBC-2} and the points are based on Eq.~\ref{eq:h(En)}.}
 \label{fig:func h(En), water}
 \end{figure}
Another interesting feature of the HSBC/MSA model is that it provides a simple approach to calculate the entropy.   Substituting $\sigma$ from Eq.~\ref{eq:NLBC-sigma} into Eq.~\ref{eq:gibbs-reac} and using $\Delta S=-\frac{\partial{\Delta G}}{\partial{T}}$,
we obtain
\begin{align} 
\Delta S=-\frac{1}{2}CRq^2 \frac{\partial{G}}{\partial{n}} \frac{\partial}{\partial{T}} \left(\frac{\hat{\epsilon}}{1+h(E_n)}\right).
\label{eq: entropy-sigma}
\end{align}
In Eq.~\ref{eq: entropy-sigma}, $\hat{\epsilon}$ and $h(E_n)$ are the only parameters that depend on temperature. We have
$\hat{\epsilon}^{\prime}=(\epsilon_{out}^{\prime}\epsilon_{out}-\epsilon_{out}^{\prime}(\epsilon_{out}-\epsilon_{in}))/(\epsilon_{out}^{2})$, where the prime indicates the derivative with respect to temperature.  The derivative of $h(E_n)$ is calculated numerically because it is a function of induced surface charge, $\sigma$, and $\sigma$ is a function of $h(E_n)$.  Finally we obtain that the HSBC/MSA solvation entropy is
\begin{align} 
\Delta S=-\frac{1}{2}CRq^2 \frac{\partial{G}}{\partial{n}} \frac{\hat{\epsilon}^{\prime}(1+h(E_n))-h^{\prime}(E_n)\hat{\epsilon}}{(1+h(E_n))^2}.
\label{eq: entropy-final}
\end{align}

The HSBC includes a minor inconsistency that deserves explanation and justification: we have used a linear-response
model to calculate the reaction potential, even though the surface charge is assumed to obey a nonlinear boundary
integral equation.  For a monovalent cation, the formally correct expression for the charging free energy is
\begin{align}
\Delta G = \int_{0}^1 \phi^{Reac}(q) dq,
\end{align}
and given the specified $h$, we have approximately
\begin{align}
\Delta G = R \int_0^1 \frac{q\;dq}{1+\frac{\delta}{R}\sqrt{|q|}}
\end{align}
where $R$ is the ion radius (the approximation is our neglect of the surface charge).  For the solvents and ions studied in this work, where $\delta \le R$, the potential still changes essentially linearly with charge.

\section{Applications}\label{sec:results}

In this section we use the HSBC/MSA model to calculate Born solvation free energies and entropies in multiple solvents, and compare our results to experimental data and the standard Born and MSA models. In Section~\ref{water}, we  present the results for water (W) and then in Section~\ref{solvents} we address methanol (MeOH), formamide (F), acetonitrile (AN), and dimethylformamide (DMF). The properties of the solvents are presented in Table~\ref{table:solvent} of the Appendix.  The temperature-dependent dielectric constant for MeOH is  
\begin{align}
  \log_{10}\epsilon&=\log_{10}\epsilon_{_{T_{0}}}-\hat{a}(T-T_0) 
  \label{eq:eps-log}
\end{align}
and the dependence for F and AN is
\begin{align}
 \epsilon&=\epsilon_{_{T_{0}}}-a(T-T_0) 
 \label{eq:eps}
 \end{align}
In both equations, $\epsilon_{T_{0}}$ is the dielectric constant of the solvent at temperature $T_0$.  The corresponding parameters are presented in Table~\ref{table:solvent}. For DMF, we used the experimental values for $\epsilon$ from Bass et. al \cite{Bass64} and fitted a cubic polynomial to the data,
\begin{align}
\epsilon_{_{DMF}}=(-1.000389\times 10^{-6})T^3+(7.718531\times 10^{-4})T^2-0.2204448T+42.04569.
\label{eq:eps-dmf}
\end{align}
\subsection{Ions in water}\label{water}

For water as a solvent, the Gibbs solvation free energies and entropies are presented in Table~\ref{results:water}. The experimental data are from Fawcett \cite{Fawcett04}, and the Born, MSA, and HSBC/MSA results are calculated. 
\begin{table}[H]
\begin{center}
\captionsetup{font=small}
\caption{Gibbs solvation free energy and entropy for ions in W at 25$^{\circ}$C, from experiment \cite{Fawcett04}, the classical Born model, the MSA Born model, and the new HSBC/MSA model. Error column represents the relative error with respect to the MSA model.}
\resizebox{\textwidth}{!}{%
\begin{tabular}{ l| c c c c c|| c c c c c}
\multirow{2}{*}{}&
\multicolumn{5}{c||}{Gibbs Energy $\Delta G$   (kJ mol$^{-1}$)}&
\multicolumn{5}{c}{Entropy $\Delta S$   (J K$^{-1}$ mol$^{-1}$)}\\
\cline{2-11}
  Ion & Expt & Born & MSA & HSBC/MSA & Error\% & Expt & Born & MSA & HSBC/MSA  & Error\%\\
  \hline
 Li$^+$ &	 $-529$ &	 $-779$ &	 $-485$ &	 $-472$ &	 $2.7$	&$-164$ 	&	 $-46$ &	 $-198$ &	 $-206$	&	 $4.0$\\

 Na$^+$ &	 $-424$ &	 $-591$ &	 $-405$ &	 $-399$ &	$1.5$ 	& $-133$ 	&	 $-35$ &	 $-142$ &	 $-148$	&	$4.2$\\

 K$^+$ &	 $-352$ &	 $-451$ &	 $-334$ &	 $-332$ &	 $0.6$ 	& $-96$ 	&	 $-27$ &	 $-100$ &	 $-103$	&	$3.0$\\

 Rb$^+$ & $-329$ &	 $-421$ &	 $-317$ &	 $-316$ &	 $0.3$ 	& $-87$ 	&	 $-25$ &	 $-91$ &	 $-94$	&	$3.3$\\

 Cs$^+$ &	 $-306$ &	 $-373$ &	 $-289$ &	 $-288$ &	 $0.3$	& $-81$ 	&	 $-22$ &	 $-77$ &	 $-79$	&	$2.6$\\
  
 F$^-$ & 	$-429$ &	 $-576$ &	 $-398$ &	 $-392$ &	 $1.5$	& $-115$ 	&	 $-34$ &	 $-138$ &	 $-143$	&	$3.6$\\
  
 Cl$^-$ &	 $-304$ &	 $-411$ &	 $-311$ &	 $-310$ &	 $0.3$	& $-53$ 	&	 $-24$ &	 $-88$ &	 $-91$	&	$3.4$\\
    
 Br$^-$ &	 $-278$ &	 $-377$ &	 $-291$ &	 $-291$ &	 $0.0$	& $-37$ 	&	 $-22$ &	 $-78$ &	 $-80$	&	$2.6$\\
      
 I$^-$ &	 $-243$ &	 $-333$ &	 $-264$ &	 $-264$ &	 $0.0$	& $-14$ 	&	 $-20$ &	 $-66$ &	 $-67$	&	$1.5$\\
\end{tabular}}
\label{results:water}
\end{center}
\end{table}
It can be observed from Table~\ref{results:water} that the HSBC/MSA results reproduce the solvation free energy of MSA model with the maximum relative error is 2.7\%, and with 4.2\% maximum relative error in the entropies.  Most of this error arises from the simple model for representing $h(E_n)$ and also use of a simple linear model for the temperature dependence of $\alpha$ (see below).

\subsection{Ions in other solvents}\label{solvents}
\begin{figure}[H]
\centering
 \includegraphics[width=25pc]{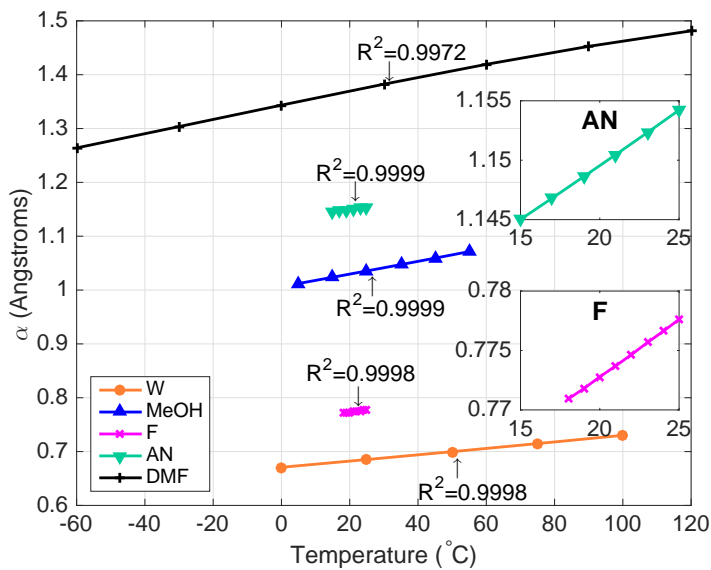}
 \caption{The linear variation of $\alpha$  with respect to temperature for different solvents. The inset plots are the enlarged graphs of $\alpha$ for acetonitrile and dimethylformamide.}
 \label{fig:alpha}
 \end{figure}

The values of $\alpha$ for methanol, formamide, acetonitrile, and dimethylformamide at different temperatures are demonstrated in Fig.~\ref{fig:alpha}. The $R^2$ values presented in this figure, clearly show that a linear relation  can describe the variation of $\alpha$ with respect to temperature. 
Writing $\alpha(T)=a_1+a_2T$, we have $h(E_n)=\left[a_1+a_2T\right]\sqrt{|E_n|}$; the parameters $a_1$ and $a_2$ for the five solvents are presented in Table~\ref{tabel:param}. 
\begin{table}[H]
\begin{center}
\caption{Values of $a_1$ and $a_2$ in the NLBC for different solvents. Values for water are reported in this table for convenience}
\begin{tabular}{ c| c c c }
Solvent	& 	$a_1$ (\AA)	&	$a_2$ (\AA~K$^{-1}$)	&	$\alpha (T=25^{\circ}$C) (\AA)	\\
\hline

W		&	0.670476		&	0.000594			& 	0.685195		\\	
MeOH	&	1.006020		&	0.001186			&	1.035502		\\	
F		&	0.753585		&	0.000960			&	0.777623		\\	
AN		&	1.131184		&	0.000920			&	1.154184		\\	
DMF		&	1.341465		&	0.001222			&	1.372015		\\	
\end{tabular}
\label{tabel:param}
\end{center}
\end{table}
The corresponding solvation free energies and entropies for MeOH, F, AN, and DMF are listed in Tables \ref{results:MeOH}--\ref{results:DMF}. The values listed as experimental results were determined using measured solvation free energies (in water) and the measured Gibbs free energies and enthalpies for transferring the ions from water to the non-aqueous solvents~\cite{Fawcett92}.
\begin{table}[H]
\begin{center}
\captionsetup{font=small}
\caption{Ion solvation free energies and entropies in MeOH at 25$^{\circ}$C. The experimental values are calculated using the Gibbs free solvation energy/entropy for water \cite{Fawcett04} and the Gibbs energy and enthalpy of transfer ions to MeOH \cite{Fawcett04}.  See Table~\ref{results:water} caption for full description. N/R: not reported.}
\resizebox{\textwidth}{!}{%
\begin{tabular}{ c| c c c c c|| c c c c c}
\multirow{2}{*}{}&
\multicolumn{5}{c||}{Gibbs Energy $\Delta G$   (kJ mol$^{-1}$) }&
\multicolumn{5}{c}{Entropy $\Delta S$   (J K$^{-1}$ mol$^{-1}$)}\\
\cline{2-11}
  Ion & Expt & Born & MSA & HSBC/MSA & Error\% & Expt & Born & MSA & HSBC/MSA & Error\%\\
  \hline
 Li$^+$ &	 $-525$ &	 $-765$ &	 $-392$ &	 $-378$ &	 $3.6 $	&	$-252$ &	 $-147$ &	 $-320$ &	 $-317$ &	 $0.9 $\\

 Na$^+$ &	 $-416$ &	 $-581$ &	 $-337$ &	 $-330$ &	 $2.1$	&	$-230$ &	 $-112$ &	 $-246$ &	 $-247$&	 $0.4 $\\

 K$^+$ & 	$-342$ &	 $-443$ &	 $-286$ &	 $-283$ &	 $1.0 $	&	$-192$ &	 $-85$ &	 $-185$ &	 $-187$&	 $1.1 $\\

 Rb$^+$ &	 $-319$ &	 $-413$ &	 $-273$ &	 $-271$ &	 $0.7 $	&	$-175$ &	 $-79$ &	 $-171$ & 	 $-173$&	 $1.2 $\\

 Cs$^+$ & $-297$ &	 $-366$ &	 $-251$ &	 $-251$ &	 $0.0 $	&	$-158$ &	 $-70$ &	 $-149$ &	 $-150$&	 $0.7 $\\
  
 F$^-$ & 	N/R  &	 $-566$ &	 $-332$ &	 $-325$ &	 $2.1 $	&	N/R &	 $-109$ &	 $-239$ &	 $-241$&	 $0.8 $\\
  
 Cl$^-$ & 	$-291$ &	 $-403$ &	 $-269$ &	 $-267$ &	 $0.7 $	&	$-69$ &	 $-78$ &	 $-166$ &	 $-168$&	 $1.2 $\\
    
 Br$^-$ & 	$-267$ &	 $-370$ &	 $-253$ &	 $-252$ &	 $0.4 $	&	$-59$ &	 $-71$ &	 $-151$ &	 $-152$&	 $0.7 $\\
      
 I$^-$ & 	$-236$ &	 $-327$ &	 $-232$ &	 $-232$ &	 $0.0 $	&	$-42$ &	 $-63$ &	 $-131$ &	 $-132$&	 $0.8 $\\ 
\end{tabular}}
\label{results:MeOH}
\end{center}
\end{table}
\begin{table}[H]
\begin{center}
\captionsetup{font=small}
\caption{Ion solvation free energies and entropies in F at 25$^{\circ}$C. The experimental values are calculated using the Gibbs free solvation energy/entropy for water \cite{Fawcett04} and the Gibbs energy and enthalpy of transfer ions to F \cite{Fawcett04}.  See Table~\ref{results:water} caption for full description. N/R: not reported.}

\resizebox{\textwidth}{!}{%
\begin{tabular}{ c| c c c c c|| c c c c c}
\multirow{2}{*}{}&
\multicolumn{5}{c||}{Gibbs Energy $\Delta G$   (kJ mol$^{-1}$) }&
\multicolumn{5}{c}{Entropy $\Delta S$   (J K$^{-1}$ mol$^{-1}$)}\\
\cline{2-11}
  Ion & Expt & Born & MSA & HSBC/MSA & Error\% & Expt & Born & MSA & HSBC/MSA & Error\%\\
  \hline
 Li$^+$ &	 $-539$ &	 $-782$ &	 $-461$ &	 $-447$ &	 $3.0$	&	$-151$ &	 $-51$ &	 $-291$ &	 $-290$ &  $0.3$\\

 Na$^+$ &	 $-432$ &	 $-593$ &	 $-388$ &	 $-382$ &	  $1.5$	&	$-162$ &	 $-39$ &	 $-210$ &	 $-211$ &  $0.5$\\

 K$^+$ &	 $-356$ &	 $-453$ &	 $-323$ &	 $-321$ &	  $0.6$	&	$-142$ &	 $-30$ &	 $-149$ &	 $-148$ &  $0.7$\\

 Rb$^+$ & $-334$ &	 $-422$ &	 $-307$ &	 $-305$ &	  $0.7$	&	$-130$ &	 $-28$ &	 $-136$ &	 $-135$ &  $0.7$\\

 Cs$^+$ &	 $-312$ &	 $-374$ &	 $-281$ &	 $-280$ &	  $0.4$	&	$-120$ &	 $-24$ &	 $-115$ &	 $-114$ &  $0.9$\\
  
 F$^-$ &	 N/R &	 $-578$ &	 $-382$ &	 $-376$ &	  $1.6$	&	N/R 	    &	 $-38$ &	 $-204$ &	 $-204$ &  $0.0$\\
  
 Cl$^-$ &	 $-290$ &	 $-412$ &	 $-302$ &	 $-300$ &	  $0.7$	&	$-87$   &	 $-27$ &	 $-131$ &	 $-130$ &  $0.8$\\
    
 Br$^-$ &	 $-267$ &	 $-378$ &	 $-283$ &	 $-282$ &	  $0.4$	&	$-78$   &	 $-25$ &	 $-117$ &	 $-115$ &  $1.7$\\
      
 I$^-$ &	 $-236$ &	 $-334$ &	 $-258$ &	 $-258$ &	  $0.0$	&	$-61$   &	 $-22$ &	 $-98$ &	 $-96$  &  $2.0$\\
\end{tabular}}
\label{results:F}
\end{center}
\end{table}
\begin{table}[H]
\begin{center}
\captionsetup{font=small}
\caption{Ion solvation free energies and entropies in AN at 25$^{\circ}$C. The experimental values are calculated using the Gibbs free solvation energy/entropy for water \cite{Fawcett04} and the Gibbs energy and enthalpy of transfer ions to AN \cite{Fawcett04}.  See Table~\ref{results:water} caption for full description. N/R: not reported.}
\resizebox{\textwidth}{!}{%
\begin{tabular}{ c| c c c c c|| c c c c c}
\multirow{2}{*}{}&
\multicolumn{5}{c||}{Gibbs Energy $\Delta G$   (kJ mol$^{-1}$) }&
\multicolumn{5}{c}{Entropy $\Delta S$   (J K$^{-1}$ mol$^{-1}$)}\\
\cline{2-11}
  Ion & Expt & Born & MSA & HSBC/MSA & Error\% & Expt & Born & MSA & HSBC/MSA & Error\%\\
  \hline
 Li$^+$ &	 $-504$ & 	$-768$ &	 $-371$ &	$-357$  &	$3.8$	&	 $-275$ &	 $-94$ &	 $-220$ &	 $-212$  &	 $3.6$ \\

 Na$^+$ &	 $-409$ & 	$-583$ &	 $-322$ &	 $-314$ &	$2.5$	&	 $-228$ &	 $-71$ &	 $-171$ &	 $-167$  &	 $2.3$\\

 K$^+$ &	 $-344$ & 	$-445$ &	 $-275$ &	 $-272$ &	$1.1$	&	 $-200$ &	 $-54$ &	 $-129$ &	 $-128$  &	 $0.8$\\

 Rb$^+$ &	 $-323$ &	 $-415$ &	 $-263$ &	 $-261$ &	$0.8$	& 	$-191$ &	 $-51$ &	 $-120$ &	 $-119$  &	 $0.8$\\

 Cs$^+$ & $-300$ &	 $-367$ &	 $-243$ &	 $-242$ &	$0.4$	& 	$-189$ &	 $-45$ &	 $-105$ &	 $-104$  &	 $1.1$\\
  
 F$^-$ & 	N/R 	    &	 $-568$ &	 $-317$ &	 $-310$ &	$2.2$	& 	N/R 	    &	 $-69$ &	 $-167$ &	 $-164$  &	 $1.8$\\
  
 Cl$^-$ & 	$-262$ &	 $-405$ &	 $-259$ &	 $-257$ &	$0.8$	& 	$-129$ &	 $-49$ &	 $-117$ &	 $-116$  &	 $0.9$\\

 Br$^-$ & 	$-247$ &	 $-371$ &	 $-245$ &	 $-244$ &	$0.4$ 	&	$-115$ &	 $-45$ &	 $-106$ &	 $-105$  & $0.9$\\
      
 I$^-$ & 	$-226$ &	 $-328$ &	 $-225$ &	 $-225$ &	$0.0$	& 	$-96$  &	 $-40$ &	 $-92$ &	 $-91$   &	 $1.1$\\ 
\end{tabular}}
\label{results:AN}
\end{center}
\end{table}
\begin{table}[H]
\begin{center}
\captionsetup{font=small}
\caption{Ion solvation free energies and entropies in DMF at 25$^{\circ}$C. The experimental values are calculated using the Gibbs free solvation energy/entropy for water \cite{Fawcett04} and the Gibbs energy and enthalpy of transfer ions to DMF \cite{Fawcett04}.  See Table~\ref{results:water} caption for full description. N/R: not reported.}
\resizebox{\textwidth}{!}{%
\begin{tabular}{ c| c c c c c|| c c c c c}
\multirow{2}{*}{}&
\multicolumn{5}{c||}{Gibbs Energy $\Delta G$   (kJ mol$^{-1}$) }&
\multicolumn{5}{c}{Entropy $\Delta S$   (J K$^{-1}$ mol$^{-1}$)}\\
\cline{2-11}
  Ion & Expt & Born & MSA & HSBC/MSA & Error\% & Expt & Born & MSA & HSBC/MSA & Error\%\\
  \hline
 Li$^+$ &	 $-539$ & 	$-768$ &	 $-335$ &	$-321$ & $4.2$		&	 $-216$ &	 $-106$ &	 $-243$ &	 $-228$  & $6.2$\\

 Na$^+$ &	 $-434$ & 	$-583$ &	 $-294$ &	 $-287$ &	$2.4$	&	 $-209$ &	 $-80$ &	 $-192$ &	 $-184$ & $4.2$\\

 K$^+$ &	 $-362$ & 	$-445$ &	 $-254$ &	 $-251$ &	$1.2$	& 	$-182$ &	 $-61$ &	 $-148$ &	 $-144$ & $2.7$\\

 Rb$^+$ &	 $-339$ &	 $-415$ &	 $-244$ &	 $-242$ &	$0.8$	& 	$-176$ &	 $-57$ &	 $-138$ &	 $-134$ & $2.9$\\

 Cs$^+$ & $-317$ &	 $-367$ &	 $-227$ &	 $-226$ &	$0.4$	& 	$-161$ &	 $-51$ &	 $-121$ &	 $-118$ & $2.5$\\
  
 F$^-$ & 	N/R	   &	 $-568$ &	 $-290$ &	 $-284$ &	$2.1$	& 	N/R 	     & $-78$ &	 $-188$ &	 $-180$ & $4.3$\\
  
 Cl$^-$ & 	$-256$ &	 $-405$ &	 $-241$ &	 $-239$ &	$0.8$	&	 $-155$ &	 $-56$ &	 $-135$ &	 $-131$ & $3.0$\\
    
 Br$^-$ & 	$-242$ &	 $-371$ &	 $-229$ &	 $-228$ &$0.4$		&	 $-156$ &	 $-51$ &	 $-123$ &	 $-120$ & $2.4$\\
      
 I$^-$ & 	$-223$ &	 $-328$ &	 $-211$ &	 $-211$ &	$0.0$	& 	$-133$ &	 $-45$ &	 $-107$ &	 $-105$ & $1.9$\\
\end{tabular}}
\label{results:DMF}
\end{center}

\end{table}

It can be seen that our HSBC/MSA model reproduces MSA results for all of the solvents tested, including being inaccurate where the MSA model is inaccurate. We note that comparing the results with experimental data, the errors are significantly larger for anions, which motivated earlier work using a charge-sign-dependent $\delta_s$ for an asymmetric MSA~\cite{Fawcett92}.  

\section{Discussion} \label{sec:discussion}

%
%

We have established that the MSA's predictions of Born ion solvation thermodynamics can be used to dramatically improve the accuracy of Poisson-based electrostatic models, as well as to enable accurate calculations of solvation entropies using dielectric theory.  Using a boundary-integral equation approach, we derived a multiscale hydration-shell boundary condition (HSBC) for the solute--solvent interface, assuming a sharp dielectric boundary.  The resulting HSBC/MSA provides solvation free energies and entropies comparable to those from the MSA model, without specific reference to solute radius.  Specifically, the HSBC/MSA model adds a simple nonlinear perturbation of the usual dielectric boundary condition,  involving the square root of the local electric field at the boundary (Eq.~\ref{eq:h(En)}). The fact that the HSBC/MSA model can reproduce the MSA results using only the normal electric field, without relying on the solute radius, suggest it may offer at least qualitative improvements for more complicated molecules such as proteins or DNA. Such applications are underway and will be reported in future work.

One of the remarkable features of the HSBC/MSA is that it furnishes a systematic approach to obtain entropies and other temperature-dependent phenomena.  A key failure of classical Poisson models in capturing temperature-dependent effects arises from the fact that the only clear parameter to vary is the solvent dielectric constant~\cite{Fawcett04}.  However, for high-dielectric solvents such as water, the $(1/\epsilon_{out}- 1/\epsilon_{in})$ factor in the Born expression is dominated by the second term, even when there are large relative changes in $\epsilon_{out}$.  The temperature dependence of $\delta_s$ has a much larger effect, which leads to semi-quantitative agreement with experimental measurements.  As a result, the continuum dielectric model using the HSBC/MSA boundary condition also accurately reproduces ion solvation entropies. Our approach offers a possible advantage over the method of Elcock and McCammon, who introduced temperature-dependent radii for complex biomolecules~\cite{Elcock97}.  More specifically, using the derivative of the boundary-integral equation Eq.~\ref{eq: entropy-final}, we can calculate entropies using a single geometry (or mesh, for numerical simulations).  In contrast, temperature-dependent radii necessitate multiple independent simulations and delicate numerical differences due to the small changes in radii.  

Using the HSBC/MSA model, we calculated the Gibbs solvation free energy and entropy for alkali metal cations and halide
anions in five solvents: water, methanol, formamide, acetonitrile, and dimethylformamide. Several interesting features
can be observed in the results. First, the HSBC/MSA model reproduces MSA free energies and entropies with high accuracy
(maximum relative error of  4.2 percent for free energy and 6.2 percent for entropy, both for solvation of Na$^+$ in
dimethylformamide).  The discrepancy arises from our simple model for the perturbation function, $h(E_n)$, and also
modeling the HSBC/MSA parameter $\alpha$ as varying linearly with temperature. Using more fitting parameters for
$h(E_n)$ and also using a more accurate interpolation (e.g. higher-order polynomial fit), the HSBC/MSA free energies and
entropies can be as accurate as the MSA results.  Second, comparing anion entropies in protic and aprotic solvents, it
can be seen that the errors in protic solvents are larger than in aprotic solvents. This could indicate that charge-sign
hydration asymmetry is more significant in protic solvents, and furthermore it may be dependent on $R_s$.  We recently
developed a semi-empirical nonlinear HSBC that captures charge-sign asymmetry~\cite{Bardhan14_asym,Bardhan15_PIERS}.
Similar to our work here, that HSBC replaced a normal-field-dependent ``radius perturbation'' with a perturbation in the
boundary condition.  However, in that work the HSBC involves an asymmetric $\tanh$, directly addressing water hydrogens'
ability to approach the solute more closely than the larger water oxygens~\cite{Mukhopadhyay14}.  Future work will
assess the similarities between our $\tanh$ HSBC and one designed to fit the charge-sign asymmetric MSA.

The present work is only a first step towards adding molecular-scale details to continuum theories; numerous applications and extensions are of interest.  For instance, to apply to biological systems in dilute electrolytes the model must be extended to allow at least Debye-Huckel treatment, if not a more sophisticated theory.  Such extensions require confronting some of the approximations inherent to modeling complex reality with a simple boundary condition. For example, HSBC/MSA modifies the boundary condition, which leads to a total surface charge different from what is predicted in the traditional macroscopic theory.  As a result, even at large distances from the ion, the electric field fails to satisfy Gauss's law.  This poses challenges for intermolecular interactions as well as extensions to dilute electrolyte solutions.  However, we have shown that this problem can be corrected by adding a renormalization charge density at a surface about one water away from the dielectric boundary, using e.g. the Stern layer~\cite{Bardhan15_PIERS}.

\section*{Acknowledgments}

This work was supported in part by the National Institute of General Medical Sciences of the National Institutes of Health under Award Number R21GM102642.  The content is solely the responsibility of the authors.

\appendices 
\section{The physical properties of the solvents}\label{appsol}
The physical properties of the solvents which are used in this research are presented in the following table.

\begin{table}[H]\
\begin{center}
\caption{Properties of the solvents and the parameters to determine the dielectric constant at different temperatures. $R_s$ is the solvent molecule radius in angstrom. $a$, $\hat{a}$, $T_0$, and $\epsilon_{_{T_{0}}}$ are parameters for Eqs.~\ref{eq:eps-log} and \ref{eq:eps}. The dielectric constant for W and DMF are third degree polynomials and are presented in Eqs.~\ref{eq:eps-water} and \ref{eq:eps-dmf} respectively. Last column shows the temperature ranges that the dielectric constant functions are valid for.}
\begin{tabular}{ l c c c c  c}

Solvent &	 $R_s $ (\AA) &	 a or $\hat{a} (10^{-2})$ &	 $T_0$ ($^\circ$C) & $\epsilon_{_{T_{0}}}$ & Range ($^\circ$C)\\
\hline

W		& 	1.420 \cite{Fawcett04} &  --- 				& --- 				& ---					& 0~--~100 \cite{Malmberg56}\\

MeOH	&	1.855 \cite{Fawcett92} & 0.26 \cite{Maryott51}	& 25	\cite{Maryott51}	& 32.63 \cite{Maryott51}	&  5~--~55 \cite{Maryott51}\\

F		& 	1.725 \cite{Fawcett92} & 72 \cite{Maryott51}	& 20 \cite{Maryott51}	& 109 \cite{Maryott51}	& 18~--~25 \cite{Maryott51}\\

AN		& 	2.135 \cite{Fawcett92} & 16 \cite{Maryott51}	& 20	\cite{Maryott51}	& 37.50 \cite{Maryott51}	& 15~--~25 \cite{Maryott51}\\

DMF		&	2.585 \cite{Fawcett92} & --- 				& --- 				& --- 				& $-60$~--~120 \cite{Bass64}\\
\end{tabular}
\label{table:solvent}
\end{center}
\end{table}

\bibliographystyle{tMPH}
\bibliography{bardhan-lab}

\begin{thebibliography}{45}
\providecommand{\url}[1]{\texttt{#1}}
\providecommand{\urlprefix}{URL }
\markboth{Taylor \& Francis and I.T. Consultant}{Molecular Physics}

\bibitem{Blum75}
L. Blum,  Mol. Phys.  \textbf{30} (5), 1529 (1975).

\bibitem{Blum78}
L. Blum and J. H{\o}ye,  Journal of Statistical Physics  \textbf{19} (4), 317
  (1978).

\bibitem{Stell75}
G. Stell and S.F. Sun,  J. Chem. Phys.  \textbf{63} (12), 5333 (1975).

\bibitem{Blum80}
L. Blum,  J. Stat. Phys.  \textbf{22} (6), 661 (1980).

\bibitem{BlumVericat80}
F. Vericat and L. Blum,  J. Stat. Phys.  \textbf{22} (5), 593 (1980).

\bibitem{Tang03}
Y. Tang,  The Journal of chemical physics  \textbf{118} (9), 4140 (2003).

\bibitem{Kirkwood34}
J.G. Kirkwood,  J. Chem. Phys.  \textbf{2}, 351 (1934).

\bibitem{Tanford57}
C. Tanford and J.G. Kirkwood,  J. Am. Chem. Soc.  \textbf{59}, 5333 (1957).

\bibitem{Sharp90}
K.A. Sharp and B. Honig,  Annu. Rev. Biophys. Bio.  \textbf{19}, 301 (1990).

\bibitem{Baker01}
N.A. Baker, D. Sept, M.J. Holst and J.A. Mc{C}ammon,  Proc. Natl. Acad. Sci.
  USA  \textbf{98}, 10037 (2001).

\bibitem{Yokota11}
R. Yokota, J.P. Bardhan, M.G. Knepley, L.A. Barba and T. Hamada,  Comput. Phys.
  Commun.  \textbf{182}, 1272 (2011).

\bibitem{Ren12_review}
P. Ren, J. Chun, D.G. Thomas, M.J. Schnieders, M. Marucho, J. Zhang and N.A.
  Baker,  Quart. Rev. Biophys.  \textbf{45}, 427 (2012).

\bibitem{Bardhan12_review}
J.P. Bardhan,  Computational Science and Discovery  \textbf{5}, 013001 (2012).

\bibitem{Borukhov97}
I. Borukhov, D. Andelman and H. Orland,  Phys. Rev. Lett.  \textbf{79}, 435
  (1997).

\bibitem{Vlachy99}
V. Vlachy,  Annu. Rev. Phys. Chem.  \textbf{50}, 145 (1999).

\bibitem{Grochowski07}
P. Grochowski and J. Trylska,  Biopolymers  \textbf{89}, 93 (2007).

\bibitem{Harris14_Boschitsch_Fenley}
R.C. Harris, A.H. Boschitsch and M.O. Fenley,  J. Chem. Phys.  \textbf{140},
  075102 (2014).

\bibitem{Gillespie15}
D. Gillespie,  Microfluid Nanofluid  \textbf{18}, 717 (2015).

\bibitem{Sandberg02_first_LD}
L. Sandberg and O. Edholm,  J. Chem. Phys.  \textbf{116}, 2936 (2002).

\bibitem{Gong09}
H. Gong and K.F. Freed,  Phys. Rev. Lett.  \textbf{102} (057603) (2009).

\bibitem{Hu12_Wei_nonlinear_Poisson_BJ}
L. Hu and G.W. Wei,  Biophys. J.  \textbf{103}, 758 (2012).

\bibitem{Abrashkin07}
A. Abrashkin, D. Andelman and H. Orland,  Phys. Rev. Lett.  \textbf{99}, 077801
  (2007).

\bibitem{Dogonadze74}
R.R. Dogonadze and A.A. Kornyshev,  J. Chem. Soc. Faraday Trans. 2
  \textbf{70}, 1121 (1974).

\bibitem{Hildebrandt04}
A. Hildebrandt, R. Blossey, S. Rjasanow, O. Kohlbacher and H.P. Lenhof,  Phys.
  Rev. Lett.  \textbf{93}, 108104 (2004).

\bibitem{Fedorov07}
M.V. Fedorov and A.A. Kornyshev,  Molecular Physics  \textbf{105}, 1 (2007).

\bibitem{Bardhan11_pka}
J.P. Bardhan,  J. Chem. Phys.  \textbf{135}, 104113 (2011).

\bibitem{Chandler72}
D. Chandler and H.C. Andersen,  The Journal of Chemical Physics  \textbf{57},
  1930 (1972).

\bibitem{Ramirez02}
R. Ramirez, R. Gebauer, M. Mareschal and D. Borgis,  Phys. Rev. E  \textbf{66},
  031206 (2002).

\bibitem{Knepley10}
M.G. Knepley, D.A. Karpeev, S. Davidovits, R.S. Eisenberg and D. Gillespie,  J.
  Chem. Phys.  \textbf{132}, 124101 (2010).

\bibitem{Borgis12}
D. Borgis, L. Gendre and R. Ramirez,  J. Phys. Chem. B  \textbf{116}, 2504
  (2012).

\bibitem{Matyushov96}
D.V. Matyushov,  Chemical physics  \textbf{211} (1), 47 (1996).

\bibitem{Vath99}
P. Vath, M.B. Zimmt, D.V. Matyushov and G.A. Voth,  The Journal of Physical
  Chemistry B  \textbf{103} (43), 9130 (1999).

\bibitem{Rocchia01}
W. Rocchia, E. Alexov and B. Honig,  J. Phys. Chem. B  \textbf{105}, 6507
  (2001).

\bibitem{Gillespie02}
D. Gillespie, W. Nonner and R.S. Eisenberg,  Journal of Physics: Condensed
  Matter  \textbf{14}, 12129 (2002).

\bibitem{Gillespie11b}
D. Gillespie,  The Journal of Physical Chemistry Letters  \textbf{2} (10), 1178
  (2011), PMID: 26295322.

\bibitem{Blum87}
L. Blum and D. Wei,  The Journal of chemical physics  \textbf{87} (1), 555
  (1987).

\bibitem{Fawcett04}
W.R. Fawcett, \emph{Liquids, solutions, and interfaces: from classical
  macroscopic descriptions to modern microscopic details}   (Oxford University
  Press, New York, 2004).

\bibitem{Fawcett92}
L. Blum and W. Fawcett,  The Journal of Physical Chemistry  \textbf{96} (1),
  408 (1992).

\bibitem{Malmberg56}
C. Malmberg and A. Maryott,  Journal of Research of the National Bureau of
  Standards  \textbf{56}, 1 (1956).

\bibitem{Bass64}
S. Bass, W. Nathan, R. Meighan and R. Cole,  The Journal of Physical Chemistry
  \textbf{68} (3), 509 (1964).

\bibitem{Elcock97}
A.H. Elcock and J.A. McCammon,  The Journal of Physical Chemistry B
  \textbf{101} (46), 9624 (1997).

\bibitem{Bardhan14_asym}
J.P. Bardhan and M.G. Knepley,  J. Chem. Phys.  \textbf{141}, 131103 (2014).

\bibitem{Bardhan15_PIERS}
J.P. Bardhan, D. Tejani, N. Wieckowski, A. Ramaswamy and M.G. Knepley,
  Progress in Electromagnetics Research Symposium (PIERS)   (2015).

\bibitem{Mukhopadhyay14}
A. Mukhopadhyay, B.H. Aguilar, I.S. Tolokh and A.V. Onufriev,  J. Chem. Theory
  Comput.  \textbf{10}, 1788 (2014).

\bibitem{Maryott51}
A. Maryott and E. Smith,  US Government Printing Office, Washington, DC
  (1951).

\end{thebibliography}

\end{document}